\documentclass[a4j,12pt]{article}
\usepackage{amsmath,amsthm,amssymb,bm}

\numberwithin{equation}{section}

\title{Conjugate Pairs of Subfactors and Entropy for Automorphisms}

\author{Marie Choda}

\date{\today}

\begin{document}
\maketitle
\centerline{Osaka Kyoiku University, 
           Asahigaoka, Kashiwara 582-8582, Japan} 
\centerline{marie@cc.osaka-kyoiku.ac.jp}

\begin{abstract} 

Based on the fact that,  for a subfactor $N$ of a II$_1$ factor $M,$ 
the first non-trivial Jones index is  2 and then 
$M$ is decomposed as the crossed product of $N$ by an outer 
action of ${\mathbb{Z}}_2,$ we study 
pairs $ \{N, uNu^* \}$  from a view point of 
entropy for two subalgebras of $M$ with a connection to the entropy for 
automorphisms, 
where the inclusion of   II$_1$ factors $ N \subset M$ 
is given as  $M$ is the crossed product of  $N$ by a finite group of outer 
automorphisms and $u$ is   a unitary in $M.$ 
\end{abstract}

keywords : {Subfactor, entropy, conjugate.}


\section{Introduction}	

For two von Neumann subalgebras $A$ and $B$ of  a finite von Neumann algebra $M,$ 
in the previous paper \cite{Ch3} we gave a modified constant $h(A | B)$ 
of the Connes-St$\o$rmer relative entropy $H(A | B)$ in \cite{CS} (cf. \cite{PP} ). 
The aim was to see the entropy for unistochastic matrices from the viewpoint of the  operator algebras 
and we   showed  among others that $h(D | uDu^*) = H(b(u)),$  
where 
$D$ is the algebra of the diagonal  matrices in the $n\times n$ 
complex  matrices $M_n(\mathbb{C}),$ and where $H(b(u))$ is  the entropy in \cite{ZSKS} 
for the unistochastic matrix $b(u)$ 
induced by a unitary $u$ in $M_n(\mathbb{C}),$ and  
in general, it does not holds  that $H(D | uDu^*) = H(b(u))$  (see, for example  \cite{PSW} ). 

In this paper, we  replace the type I$_n$ factor $M_n(\mathbb{C})$ to a II$_1$ factor $M$. 
The above relation in \cite{Ch3} (cf. \cite{Oka}) suggests us 
that if $A$ and $B$ are maximal abelian subalgebras 
of a II$_1$ factor $M$, then $h(A|B)$ is not necessarily finite.  
In order  to discuss  on a subalgebra  $A$  of $M$ 
with   $h(A|uAu^*) < \infty$ for all unitaries $u$ in $M$,  
we pick up here a subfactor $N \subset M$ with the Jones index $[M : N] < \infty$ (\cite{J}) 
 (cf. \cite{EK}). 
\smallskip

Let $M$ be a type II$_1$ factor, and let $N$ be a subfactor of $M$ such that $[M : N] = 2$  
which is the simplest, nontrivial unique subfactor of $M$ up to conjugacy. 
Then $M$ is decomposed into 
the crossed product of $N$ by an outer automorphism with the period 2. 
Based on this fact, 
we study  the set of the values $h(N | uNu^*)$ for the inclusion of factor-subfactor $N \subset M,$ 
with a connection to the inner automorphisms $Ad u,$ 
where $M$ is given as the crossed product  $N \rtimes_\alpha G$ 
of a type II$_1$ factor $N$ by a finite group $G$  with respect to an outer action $\alpha$ 
and $u$ is a unitary in $M.$ 

First, for two von Neumann subalgebras $A$ and $B$ of a finite von Neumann subalgebra $M$, 
we show, in Corollary 2.2.3 below, 
that if $E_A E_B = E_B E_A$ (which is called  the commuting square  condition 
in the sense of \cite{GHJ}) then $H(A |B) = h(A | B),$ 
where $E_A$ is the conditional expectation of $M$ onto $A.$  

We give an extended notion $H_N(Ad u)$ of $H(b(u))$ to  the inner automorphisms $Ad u$ 
in Definition 3.1.1, and   
we show that $h(N | uNu^*)  \leq H_N(Ad u)$ in Theorem 3.1.4. 

The inner conjugacy class of $N$ is rich from the view point of the values of 
$h(N | uNu^*),$ that is, 
in the special  case of $G = {\mathbb{Z}}_2$, 
keeping the fact 
that $h(N \ |  uNu^*) \leq   H(M |uNu^*) = \log 2$ for all unitary $u \in M$ 
in mind, we have 
that $\{ h(N \ |  uNu^*) : u \in M \ \text{a unitary} \} = [0, \log 2]$ in 
Theorem 3.2.3. 

\vskip 0.3cm

\section{Preliminaries}
In this section, we summarize, for future reference, notations, terminologies and 
basic facts.  

Let $M$ be a finite von Neumann algebra, and let $\tau$ be a 
faithful normal tracial state. 
For each von Neumann subalgebra $A$, there is a unique $\tau$ preserving conditional 
expectation $E_A : M \to A.$ 

\subsection{\bf Connes-St\o rmer relative   entropy } 

Let $S$ be the set of all finite families $(x_i)_i$ of positive 
elements in $M$ with $1 = \sum_i x_i.$ 
Let $A$ and $B$ be two von Neumann subalgebras of $M.$ 
The  relative  entropy $H (A | B)$   is defined by Connes and St$\o$rmer(\cite{CS}) as 
$$H( A \mid B) = 
\sup_{(x_i) \in S} \sum_{i} (\tau \eta E_{B}(x_i) - \tau \eta E_A(x_i)). $$
Here,  $\eta$ is the function defined by 
$$\eta(t) = -t \log t, \ (0 < t \leq 1) \quad \text{and} \ \eta(0) = 0.$$ 
Let $\phi$ be a normal state on $M.$ 
Let $\Phi$ be the set of all finite families $(\phi_i)_{i}$ of positive 
linear functionals on $M$ with $\phi = \sum_i \phi_i.$ 
The relative  entropy $H_\phi(A | B)$ of $A$ and $B$ with respect to $\phi$ 
is given by Connes (\cite{Co2}) as 
$$ H_\phi (A | B) = 
\sup_{(\phi_i) \in \Phi} 
\sum_i ( S(\phi_i\mid_A, \phi \mid_A) - S(\phi_i\mid_B, \phi \mid_B) )$$ 
and  if $\phi = \tau$ then 
$H_\tau(A | B) = H( A \mid B).$ 
Here $S(\phi, \psi)$ is the relative entropy for  positive linear functionals  $\phi$ and $\psi$ 
on $M$ (cf. \cite{NS, OP}). 
\smallskip

\subsection{Modified relative entropy  for two subalgebras. }

We modified in \cite{Ch3} the Connes-St\o rmer relative entropy for a pair of 
subalgebras  as follows : 

Let $A$ and $B$ be two von Neumann subalgebras of $M.$ 
Let $S$ be the set of all finite families $(x_i)$ of positive 
elements in $M$ with $1 = \sum_i x_i.$ 
The {\it conditional relative entropy} $ h( A \mid B)$ 
of  $A$ and $B$  corresponding $H(A | B)$ is given as 
$$ h( A \mid B) 
= \sup_{(x_i) \in S}  
\sum_{i} (\tau \eta E_{B}(E_{A}(x_i)) - \tau \eta E_{A} (x_i)). $$ 

Let $S(A) \subset S$ be the set of all finite families $(x_i)$ of positive 
elements in $A$ with $1 = \sum_i x_i.$ 
Then it is clear that 
$$ h( A \mid B) 
= \sup_{(x_i) \in S(A)}  
\sum_{i} (\tau \eta E_{B}(x_i) - \tau \eta (x_i)). $$
\vskip 0.3cm

Let $S'(A) \subset S(A)$ be the set of all finite families $(x'_i)$ 
with each $x'_i$ a scalar multiple of a projection in $A.$ 
Then 
$$ h( A \mid B) 
= \sup_{(x'_i) \in S'(A)}  
\sum_{i} (\tau \eta E_{B}(x'_i) - \tau \eta (x'_i)). $$ 

The {\it conditional relative entropy} of 
$A$ and $B$ with respect to $\phi$ corresponding $H_\phi(A | B)$ is given as 
$$ h_\phi(A|B) = 
\sup_{(\phi_i) \in \Phi} 
 \sum_i ( S(\phi_i\mid_A, \phi \mid_A) - 
S((\phi_i\circ E_A)\mid_B, (\phi \circ E_A)\mid_B ).$$ 
\smallskip
If we let $\Phi(A) \subset \Phi$ be the set of all finite families $(\phi_i')$ 
in $\Phi$ with each $\phi_i' = \phi_i' \circ E_A,$ then 
$$ h_\phi(A|B) = 
\sup_{(\phi_i') \in \Phi(A)} 
 \sum_i ( S(\phi_i'\mid_A, \phi \mid_A) - 
S(\phi_i' \mid_B, (\phi \circ E_A)\mid_B ).$$ 
\smallskip

We give conditions for that  $h_\phi(A | B) = H_\phi(A | B)$ in Corollary 2.2.4 below, 
and  show relations for $h_\phi(A | B), H_\phi(A | B), H_\phi(A)$ and $h_\phi(A)$ in Theorem 2.2.2, 
where $h_\phi(A)$ is given   by modifing 
$H_\phi(A)$ in \cite{Co2} for a von Neumann subalgebra $A$ of $M$ (cf. \cite{NS}), 
that is 
$$h_\phi(A) = \sup_{(\phi_i) \in \Phi(A)}  
\sum_i ( \eta(\phi_i(1)) + S(\phi_i|_A, \phi|_A) ) $$ 
and 
$$H_\phi(A) = \sup_{(\phi_i) \in \Phi}  
\sum_i ( \eta(\phi_i(1)) + S(\phi_i|_A, \phi|_A) ). $$ 

Cleary, we have that 
$0 \leq h_\phi(A) \leq H_\phi(A).$  
In the case of $\phi$ is the  trace $\tau,$ 
$h_\phi(A) = h(A),$ which is given as 
$$h(A) = \sup_{(x_i) \in S(A)}  
\sum_i ( \eta(\tau(x_i)) - \tau \eta(x_i) ). $$ 
\smallskip

We need the following lemma in order to prove  Theorem 2.2.2, in which we show relations 
among  $H_\phi(A), h_\phi(A), H_\phi(A | B)$ and  $ H_\phi(A | B) $ : 
\vskip 0.3cm

\noindent{\bf Lemma 2.2.1. } 
{\it Let $A$ and $B$ be von Neumann subalgebras of  a finite von Neumann algebra $M$, 
and let $\psi, \phi$ be 
positive linear functionals on $M.$ 
If $E_A E_B = E_B E_A,$ then 
$$S((\psi \circ E_A) \mid_B, (\phi \circ E_A) \mid_B) = 
S((\psi \circ E_A) \mid_{A \cap B}, (\phi \circ E_A) \mid_{A \cap B}).$$ }

 \begin{proof} 
Relative entropy for positive linear functionals $\psi, \phi$ on a unital 
$C^*$-algebra $C$ is expressed in \cite {K} by 
$$S(\psi, \phi) = 
\sup_{n \in \mathbb{N}} \sup_x \{\psi(1)\log n - \int_{1/n}^\infty 
\left( \psi (y(t)^*y(t)) + 
\frac 1t \phi (x(t)x(t)^*) \right) \frac{dt }t\}$$
where $x(t) : (\frac 1n, \infty) \to C$ is a step function with finite range,  
and $y(t) = 1 - x(t).$ 

Let  $x(t) : (\frac 1n, \infty) \to B$ be a  step function with finite range, 
then 
$E_A(x(t)) : (\frac 1n, \infty) \to A \cap B$ is a step function 
with finite range because of that $E_A E_B = E_B E_A.$ 
Since $E_A(x)^* E_A(x) \leq E_A(x^*x)$ for all $x \in M,$ 
we have that 
\begin{eqnarray*}
\lefteqn{
\psi(1)\log n - \int_{1/n}^\infty 
\left( \psi \circ E_A(y(t)^*y(t)) + 
\frac 1t \phi \circ E_A(x(t)x(t)^*) \right) \frac{dt }t 
}\\
&&\leq \psi(1)\log n - \int_{1/n}^\infty 
\left( \psi (E_A(y(t))^* E_A(y(t))) + 
\frac 1t \phi (E_A(x(t)) E_A(x(t)^*)) \right) 
\frac{dt }t \\
&& \leq S((\psi \circ E_A) \mid_{A \cap B}, 
(\phi \circ E_A) \mid_{A \cap B}).
\end{eqnarray*} 
This implies that 
$$S((\psi \circ E_A) \mid_B, (\phi \circ E_A) \mid_B) \leq
S((\psi \circ E_A) \mid_{A \cap B}, (\phi \circ E_A) \mid_{A \cap B}).$$ 
Since the opposite inequality is clear, 
we have the equality. 
\end{proof} 
\vskip 0.3cm
\smallskip

\noindent{\bf Theorem 2.2.2.} 
{\it Let $M$ be a finite von Neumann algebra with a  normal faithful tracial state $\tau.$ 
Let $\phi$ be a normal state of  $M, and $ 
let $A, B$ be von Neumann subalgebras  of $M.$ 
Then }
\smallskip

\noindent (1) 
$h_\phi(A | B) \le  h_\phi(A |{\mathbb{C}}{\bf 1}) = h_\phi(A).$ 
\smallskip

\noindent (2)  
{\it Assume that $E_A E_B = E_B E_A.$ Then 
$h_\phi(A | B) =  h_\phi(A | A \cap B). $

Hence,  if  $A \cap B = \mathbb{C},$ then 
$h_\phi(A | B) = h_\phi(A).$} 
\smallskip

\noindent (3) 
{\it If  $E_A E_B = E_B E_A$ and if $\phi = \phi \circ E_A,$ then 
$H_\phi(A | B) =  H_\phi(A | A \cap B).$ 

Especially,  if  $A \cap B = \mathbb{C},$ then 
$H_\phi(A | B) = H_\phi(A).$ }
\smallskip

\noindent (4) 
{\it 
If $B \subset A,$ then 
$$H_\phi(A | B) = h_\phi(A | B).$$ 

Especially, $H_\phi(A) = h_\phi(A).$}
\vskip 0.3cm

\begin{proof} 
(1) Let $\phi, \psi$ be normal states of $M,$ and let $(\phi_i)_i \in \Phi.$ 
Then $(\phi \circ E_A) \mid_B$ and $\frac 1 {\phi_i (1)} (\phi_i \circ E_A)\mid_B$ 
are  states of $B$ so that 
$$S(\frac 1 {\phi_i (1)} (\phi_i \circ E_A)\mid_B, \ (\phi \circ E_A)\mid_B) 
\geq 0$$
and 
\begin{eqnarray*}
\lefteqn{S(\frac{ 1}{\phi_i (1)} (\phi_i \circ E_A)\mid_B, \ (\phi \circ E_A)\mid_B)} \\
& & = \frac 1{ \phi_i (1)} S(\phi_i \circ E_A \mid_B, \ (\phi \circ E_A)\mid_B ) - \phi_i (1) \eta(\frac 1 { \phi_i (1)} ) \\
& & = \frac 1{ \phi_i (1)} S(\phi_i \circ E_A \mid_B, \ (\phi \circ E_A)\mid_B) -  \log( \phi_i (1)). 
\end{eqnarray*} 
Hence 
$$-S((\phi_i  \circ E_A) \mid_B, \ (\phi \circ E_A)\mid_B ) 
\leq \eta(\phi_i(1)),$$ 
and if $B = \mathbb{C}{\bf 1}$ then the equality holds 
because 
$\frac 1{\phi_i(1)} \phi_i|_{\mathbb{C}{\bf 1}} 
= \phi|_{\mathbb{C}{\bf 1}}. $ 
These imply that 
\begin{eqnarray*}
\lefteqn{h_\phi(A|B)}\\ 
&& = \sup_{(\phi_i) \in \Phi}  
            \sum_i (S(\phi_i\mid_A, \phi \mid_A) 
           - S((\phi_i\circ E_A)\mid_B, (\phi \circ E_A)\mid_B)) \\
 & & \leq \sup_{(\phi_i) \in \Phi } 
       \sum_i (S(\phi_i\mid_A, \ \phi \mid_A) + \eta(\phi_i(1))) \\
&& = h_\phi(A),
\end{eqnarray*} 
and the equality holds if $B = \mathbb{C}{\bf 1}.$ 
\vskip 0.3cm

(2) Assume that $E_A E_B = E_B E_A$, then by  lemma 2.2.1, we have that 
\begin{eqnarray*}
\lefteqn{
h_\phi(A|B) 
}\\
&& = \sup_{(\phi_i) \in \Phi}  \sum_i ( S(\phi_i\mid_A, \phi \mid_A) 
- S((\phi_i\circ E_A)\mid_B, \ (\phi \circ E_A)\mid_B))\\
&& = \sup_{(\phi_i) \in \Phi}  \sum_i ( S(\phi_i\mid_A, \phi \mid_A) 
- S((\phi_i\circ E_A)\mid_{A \cap B}, \ (\phi \circ E_A)\mid_{A \cap B}))\\
&& = h_\phi(A| A \cap B) . 
\end{eqnarray*}
Since $h_\phi(A|B)$ is decreasing in $B,$ it implies that 
$$h_\phi(A|B) = h_\phi(A| A \cap B).$$ 
Hence, if  $A \cap B = \mathbb{C}$ 1,  then 
$h_\phi(A | B) = h_\phi(A| \mathbb{C} 1) = h_\phi(A).$ 
\smallskip 

(3) 
Assume that $E_A E_B = E_B E_A$ and that $\phi \circ E_A = \phi.$ 
Let $(\phi_i)_i \in \Phi,$ then $(\phi_i \circ E_A)_i  \in \Phi$ and we have that 
\begin{eqnarray*}
\lefteqn{
H_\phi(A|B) 
}\\
&& = \sup_{(\phi_i) \in \Phi}  \sum_i ( S(\phi_i\mid_A, \phi \mid_A) 
- S(\phi_i \mid_B, \ \phi \mid_B))\\
&&\geq \sup_{(\phi_i) \in \Phi}  
    \sum_i ( S((\phi_i\circ E_A)\mid_A, \phi \mid_A) 
    - S((\phi_i\circ E_A)\mid_B, \ \phi \mid_B))\\
&& = \sup_{(\phi_i) \in \Phi}  \sum_i ( S(\phi_i\mid_A, \phi \mid_A) 
- S((\phi_i\circ E_A)\mid_B, \ (\phi \circ E_A)\mid_B))\\
&& = \sup_{(\phi_i) \in \Phi}  \sum_i ( S(\phi_i\mid_A, \phi \mid_A) 
- S((\phi_i\circ E_A)\mid_{A \cap B}, \ (\phi \circ E_A)\mid_{A \cap B}))\\
&& = \sup_{(\phi_i) \in \Phi}  \sum_i ( S(\phi_i\mid_A, \phi \mid_A) 
- S(\phi_i \mid_{A \cap B}, \ \phi \mid_{A \cap B}))\\
&& = H_\phi(A| A \cap B) . 
\end{eqnarray*}
In general, $H_\phi(A| A \cap B) \geq H_\phi(A| B) $ so that 
$$H_\phi(A| B) = H_\phi(A| A \cap B).$$ 
Especially,  $H(A| B) = H(A| A \cap B)$ (cf. \cite{WW}), and 
 if  $A \cap B = \mathbb{C},$ then 
$H_\phi(A| A \cap B) = H_\phi(A| {\mathbb{C}} 1 ) = H_\phi(A)$ 
so that 
$$H_\phi(A | B) = H_\phi(A) .$$ 
\smallskip 

(4) If $B \subset A,$ then 
\begin{eqnarray*}
\lefteqn{
H_\phi(A|B) 
}\\
&& = \sup_{(\phi_i) \in \Phi}  \sum_i ( S(\phi_i\mid_A, \phi \mid_A) 
- S(\phi_i \mid_B, \ \phi \mid_B))\\
&& = \sup_{(\phi_i) \in \Phi}  \sum_i ( S(\phi_i\mid_A, \phi \mid_A) 
- S((\phi_i\circ E_A)\mid_B, \ (\phi \circ E_A)\mid_B))\\
&& = h_\phi(A| B) . 
\end{eqnarray*}
By combining with (1), we have $H_\phi(A) = h_\phi(A).$ 
\end{proof} 
\vskip 0.3cm

\noindent{\bf Corollary 2.2.3.} 
{\it Assume that $E_A E_B = E_B E_A.$ Then 
$H(A | B) = h(A | B)$. 
Moreover, if  $\phi = \phi \circ E_A,$ then 
$$H_\phi(A | B) = h_\phi(A | B).$$}

\begin{proof} Assume that $E_A E_B = E_B E_A$ and  that $\phi = \phi \circ E_A.$ 
Then 
by using Theorem 2.2.2 (3), (4) and (2), we have that 
$$H_\phi(A | B) = H_\phi(A | A \cap B) = h_\phi(A | A \cap B)
  = h_\phi(A | B).$$
Since  $E_A$  is the $\tau$-conditional expectation, $\tau = \tau \circ E_A,$ hence 
$$H(A | B) = h(A | B).$$
\end{proof} 

\vskip 0.3cm
\section{Inner conjugate subfactors } 

Connes-St$\o$rmer defined  
the entropy $H(\alpha)$  for a trace preserving automorphism 
$\alpha$ of a finite von Neumann algebra in \cite{CS}. The definition is arivable for 
a trace preserving *-endomorphism too. 

For a trace preserving *-endomorphism $\sigma$ of a finite von Neumann 
algebra $N,$ it was shown 
a relation between the entropy $H(\sigma)$ for  $\sigma$ 
and the relative entropy $H(N \mid \sigma(N))$  in the papers 
\cite  {Ch1, Ch2, Ch-H, H, St} (cf. \cite{NS}). 
The relation is, roughly speaking, that 
$$H(\sigma) = \frac 12 H(N \mid \sigma(N))$$ 
under a certain condition. 
Such a *-endomorphism $\sigma$ can be extended offten to an 
automorphism $\alpha$ of a finite von Neumann algebra $M$ which contains  
$N$ as a von Neumann subalgebra. 
Some examples of such endomorphisms appeared in a relation to Jones index 
theory of subfactors. 
In \cite{Ch2}, we studied a nice class of such a *-endomorphism 
$\sigma$ of a type II$_1$factor $N$ which is extendable to an automorphism $\alpha$ 
of the big type II$_1$factor $M$ obtained by the basic construction from 
$N \supset \sigma(N).$ 
We called such a $\sigma$ {\it basic} *-endomorphism and showed that 
$H(\alpha) = \frac 12 H(N | \sigma(N)).$ 
Since $\sigma(N) \subset N,$  we have by Theorem 2.2.2 (4) that $H(N | \sigma(N)) = h(N | \sigma(N))$ 
so that 
$$H(\alpha) = \frac 12 H(N | \sigma(N)) = \frac 12 h(N | \sigma(N)) = \frac 12 h(N |\alpha(N)).$$
This means that for an automorphism $\alpha$ of a II$_1$ factor 
$M$ we may be able to choose a subfactor $N \subset M$ such that 
the entropy for $\alpha$  is given from 
$h(N |\alpha(N)).$ 

Our study in this section is motivated by these results. 
The  above automorphism $\alpha$ arising from a *-endomorphism as is outer. Here,  
we discuss by replacing the $\alpha$ to inner  automorphisms $Ad u$ and the entropy $H(\alpha)$ to 
the entropy $H_N(Ad u)$ defined below. 
\vskip 0.3cm

\subsection{Entropy for Inner Automorphisms with respect to Subfactors} 
Let $N$ be a type II$_1$ factor with the canonical trace $\tau$ and let 
$G$ be a finite group. Let $\alpha$ be an outer action of $G$ on $N,$ 
so that for all $g \in G, g \ne 1$ 
if $\alpha_g(x)a = ax$ for all $x \in N,$ then $a = 0.$ 
Hereafter, we let $M$  be the crossed product of $N$ by $G$ 
with respect to  $\alpha: $  
$$ M = N \rtimes_\alpha G.$$ 
We identify $N$ with the von Neumann subalgebra embedded in $M,$ and 
denote by $v$ the  unitary representation of $G$ in $M$ 
such that  every $v_g$ is a unitary in $M$ with 
$$\alpha_g (x) = v_g x v_g^*, \quad (x \in N, g \in G).$$ 
Then every  $x \in M$ is written by the Fourier expansion 
$$x = \sum_{g \in G} x_g v_g, \quad (x_g \in N) $$ 
and $x_g = E_N(xv_g^*).$ 
A $u \in M$ is a unitary if and only if  
$$\sum_{g \in G} u_{hg} \alpha_h(u_g^*) = \delta_{h, 1} 
\quad \text{and} \quad  
\sum_{g \in G} \alpha_{g}^{-1}(u_g^* u_{gh}) = \delta_{h, 1}, $$
where we denote the identity of $G$ by $1.$ 
This imply that $\sum_{g \in G} \tau(u_gu_g^*) = 1,$ 
and we can put as the followings : 
\vskip0.3cm

\noindent{\bf Definition 3.1.1.} 
The {\it entropy of the inner automorphism $\rm{Ad} u$ of $M$ 
with respect to $N$}  is given by 
$$H_N(Adu) = \sum_{g \in G} \eta\tau(u_gu_g^*).$$ 
\vskip0.3cm

\noindent{\bf Comment 3.1.2.} 
Each $x \in M$ is representated as the matrix  $x = (x(g,h))_{gh}$ indexed by the elements of $G.$ 
Here $x(g,h) \in N$ for all $g,h$ in $G,$ and 
$x(g,h) = \alpha_g^{-1} (E_N(x v_h^*)) = \alpha_g^{-1}(x_h).$ 
The entropy $H(b(u))$ defined in \cite{ZSKS} is written as 
$$H(b(u)) = \frac 1n \sum_{i,j} \eta( |u(i, j)|^2 ),$$ 
when $b(u)$ is the unistochastic matrix induced by a unitary 
$u = (b(i,j))_{ij}$  in  $M_n(\mathbb{C}).$ 
A matrix representation for an $x$  in  $M_n(\mathbb{C})$ is depend on the diagonal matrix 
algebra. 
In that sense, we consider the notion of $H_N(Adu)$ corresponds to the  notion of the entropy 
for a unistochastic matrix. 
\vskip0.3cm

\noindent {\bf Lemma 3.1.3.} 
\noindent (1) {\it If  Ad$u$ and Ad$w$ are conjugate, then 
$H_N(Adu) = H_N(Adw).$}
\smallskip 

\noindent (2) {\it 
If $\theta = Adu$ for some unitary $u \in M,$ 
then 
$H_N(\theta^{-1} ) = H_N(\theta).$ }

\begin{proof}
(1) Assume $Adu = \theta Adw \theta^{-1}$ for some automorphism $\theta$ of $M.$ 
Then $\theta(w) = \lambda u$ for some complex number $\lambda$ with 
$|\lambda| = 1$ and so 
$\eta \tau(w_gw_g^*) = \eta \tau(u_gu_g^*)$ 
which implies that $H_N(Adu) = H_N(Adw).$ 

(2) 
Let $w \in M$ be a unitary with $Adw = \theta^{-1},$ then $w = \gamma u^*$ 
for some $\gamma \in \mathbb{T}.$ For the expression that 
$w = \sum_g w_g v_g,$ we have that $w_g = \gamma \alpha_g(u_{g^{-1}}^*)$ 
for all $g \in G$ so that 
$$H_N(\theta^{-1} ) = 
\sum_g \eta\tau(w_gw_g^*) =  \sum_g \eta\tau(u_gu_g^*) = H_N(\theta).$$
\end{proof}
\smallskip

The  $h(N | uNu^*)$ is bounded by $ H_N(Adu)$  as follows : 
\smallskip

\noindent {\bf Theorem 3.1.4.} 
{\it Assume that
$N$ is a type II$_1$ factor, $G$ is  a finite group 
and $M = N \rtimes_\alpha G$ with respect to  the outer action $\alpha.$ 
Then  for each unitary $u \in M,$ we have that 
$$h(N | uNu^*) \leq H_N(Adu) = \sum_{g \in G} \eta \tau(u_gu_g^*) \leq \log |G|,$$ }
where $|G|$ is the cardinarity of $G.$ 
\begin{proof}
Let $(\lambda_i p_i)_{i \in I} \in S'(N)$ be a finite partition of 
the unity, that is, 
$$\sum_{i \in I} \lambda_i p_i = 1$$
where $(\lambda_i)_{i \in I} $ are positive numbers and 
$(p_i)_{i \in I} $ are projections in $N.$ 
For a given $\varepsilon,$ choose an $\epsilon > 0$ with 
with $2|G|\eta(\epsilon) < \min\{\varepsilon, 1/e\}.$  
There exist mutually orthogonal projections $(q_{i,k})_k \subset N$ 
and nonnegative numbers $(\alpha^g_{i,k})_k$ 
which satisfy that 
$$p_i = \sum_k q_{i,k} \quad \text{and }
\quad 0 \leq q_{i,k} u_gu_g^* q_{i,k} - \alpha^g_{i,k} q_{i,k} \leq 
\epsilon q_{i,k}.$$ 
This is possible by the induction method 
of the spectral decompositions 
for $(p_iu_gu_g^*p_i)_{i \in I, g \in G},$  
(see for example. \cite [Proof of 4.3 Lemma] {PP}). 
In fact, letting $G = \{g_1, \cdots, g_m\}$ and 
by the spectral decomposition for $p_iu_{g_1}u_{g_1}^*p_i \in p_iNp_i,$ 
we have mutually orthogonal projections $(q_{i,k}^1)_{k_1} \subset p_iNp_i$ 
and nonnegative numbers $(\alpha^1_{i,k_1})_{k_1}$ 
$$p_i = \sum_{k_1} q_{i,{k_1}}^1 \quad \text{and }
\quad 0 \leq 
q_{i,{k_1}}^1u_{g_1}u_{g_1}^* q_{i,{k_1}}^1 - 
\alpha^1_{i,{k_1}} q_{i,{k_1}}^1 \leq 
\epsilon q_{i,{k_1}}^1.$$ 

Next by the consideration  for  $q_{i,k_1}^1u_{g_2}u_{g_2}^* q_{i,k_1}^1$,  
we have a partition $(q_{i,k_1,k_2}^2)_{k_2}$ of $q_{i,k_1}^1$ 
 and $(\alpha^2_{i,k_1,k_2})_{k_2}$. 
Put $\alpha^{g_j}_{i,k} = \alpha^j_{i,k_1, \cdots, k_j}$ 
and $q_{i,{k}} = q^m_{i,k_1, \cdots, k_m}.$ 
Then these satisfy the desired conditions. 

Since $\eta(x + y) \le \eta(x) + \eta(y), (x, y \in N), $ 
$\eta$ is increasing on $[0, 1/e],$ 
and the family $(q_{i,k})_k$ is mutually orthogonal,
we have for $\epsilon$ with $\epsilon \le 1/e$   
\begin{eqnarray*} 
\lefteqn{\tau \eta(\sum_k q_{i,k}u_gu_g^*q_{i,k})} \\ 
&&\leq \tau \eta(\sum_k (q_{i,k} u_gu_g^* q_{i,k} - \alpha^g_{i,k} q_{i,k})) 
 + \tau \eta(\sum_k \alpha^g_{i,k} q_{i,k})\\
&&\leq \eta(\epsilon) \tau(p_i) + \sum_k 
      \eta( \alpha^g_{i,k}) \tau( q_{i,k}).
\end{eqnarray*}
Hence, by the condition that $\sum_{i,k} \lambda_i \tau( q_{i,k}) = 1,$ 
the operator concavity of $\eta$  implies that 
\begin{eqnarray*} 
\lefteqn
{\sum_{i,g} \lambda_i \tau \eta(\sum_k q_{i,k}u_gu_g^*q_{i,k}) }\\
&&\leq  \sum_{i,k, g} \lambda_i  \eta( \alpha^g_{i,k}) \tau( q_{i,k}) 
  + |G| \eta(\epsilon)\\
&&= \sum_{i,k} \lambda_i \tau( q_{i,k}) \sum_g \eta(\alpha^g_{i,k})  
  + |G| \eta(\epsilon)\\
&&\leq \sum_g \eta(\sum_{i,k} \lambda_i \tau( q_{i,k})\alpha^g_{i,k})   
+ |G| \eta(\epsilon).
\end{eqnarray*}

Remark that for all $g \in G$ 
$$\tau(u_gu_g^*) - \sum_{i,k} \lambda_i \alpha^g_{i,k}\tau( q_{i,k}) 
= \sum_{i,k} \lambda_i \tau( q_{i,k} u_gu_g^* q_{i,k} 
- \alpha^g_{i,k} q_{i,k}) $$
and that
$$0 \leq \sum_{i,k} \lambda_i \tau( q_{i,k} u_gu_g^* q_{i,k} 
- \alpha^g_{i,k} q_{i,k}) \leq \epsilon.$$
Then we have that 
$$0 \leq \tau(u_gu_g^*) - \sum_{i,k} \lambda_i \tau( q_{i,k})\alpha^g_{i,k} 
\leq \epsilon, \quad (g \in G)$$
so that 
$$\sum_g \eta(\sum_{i,k} \lambda_i \tau( q_{i,k})\alpha^g_{i,k})  
\leq \sum_{g} \eta\tau(u_gu_g^*) + |G|\eta(\epsilon).$$
Here we used the following inequality in \cite [(2.8)] {NS} 
$$|\eta(s) - \eta(t)| \leq \eta(s-t) \ \text{for} \ 0 \leq s-t \leq \frac12.$$
Remark that  $\sum_g u_{g} u_g^* =1$ and that for all $i$ 
the projections  $(q_{i,k})_k$ is mutually orthogonal.  
Hence by using  the following fact that 
$$\tau \eta(u_g^*q_{i,k}u_g) = \tau \eta(q_{i,k}u_gu_g^*q_{i,k}), $$ 

we have that 
\begin{eqnarray*} 
\lefteqn{\sum_{i,k} 
(\tau \eta E_N(u^* \lambda_i q_{i,k} u) 
  -  \tau \eta(\lambda_i q_{i,k}))} \\
&&= \sum_{i,k} 
(\tau \eta (\sum_{g \in G} \alpha_g (u_g^* \lambda_i q_{i,k} u_g) 
  -   \eta(\lambda_i) \tau(q_{i,k}) ) )\\
&&\leq  \sum_{i,k} \sum_g 
(\tau \eta (u_g^* \lambda_i q_{i,k} u_g) 
  -  \eta(\lambda_i) \tau(q_{i,k}) ) \\
&& = \sum_{i,k,g } (\eta(\lambda_i) \tau(u_g^*q_{i,k}u_g) 
 + \lambda_i \tau \eta(u_g^*q_{i,k}u_g) 
 - \sum_{i,k}\eta(\lambda_i) \tau(q_{i,k}) \\
&& = \sum_i (\eta(\lambda_i) \tau(p_i \sum_g u_gu_g^*) 
 + \sum_{i,k,g } \lambda_i \tau \eta(u_g^*q_{i,k}u_g) 
 - \sum_{i}\eta(\lambda_i) \tau(p_i) \\
&& = \sum_{i,k,g } \lambda_i \tau \eta(q_{i,k}u_gu_g^*q_{i,k}) \\
&&= \sum_{i,g} \lambda_i \tau \eta(\sum_k q_{i,k}u_gu_g^*q_{i,k}) \\
&&\leq \sum_{g} \eta\tau(u_gu_g^*) + 2|G|\eta(\epsilon). 
\end{eqnarray*}
Thus 
\begin{eqnarray*} 
\lefteqn{h(N | u N u^*) 
= \sup_{(\lambda_i p_i)\in S'(N)} \sum_i (\tau \eta E_N(u^* \lambda_i p_i u) 
  -  \tau \eta(\lambda_i p_i) )} \\
&&= \sup_{(\lambda_i q_{ik})_{i,k}} \sum_{i,k} 
(\tau \eta E_N(u^* \lambda_i q_{i,k} u) 
  -  \tau \eta(\lambda_i q_{i,k}))\\
&&\leq \sum_g \eta \tau(u_gu_g^*).
\end{eqnarray*}
\end{proof} 
\vskip 0.3cm

Since $\eta$ is a concave function, Theorem 3.1.4 implies the following : 
\smallskip

\noindent {\bf Corollary  3.1.5.} 
{\it Assume that $N, G, u$ be as in Theorem 3.1.4 and that 
$h(N | uNu^*) = \log |G|.$ 
Then $\tau(u_gu_g^*) = \frac 1{|G|}$ for all $g \in G.$}
\vskip 0.3cm

\noindent {\bf Remark and Example 3.1.6.   } 
Let $A$ and $B$ be  subalgebras of a type II$_1$ factor $M.$ Then 
$h(A|B) \leq H(A | B) \leq H(M|B),$ 
and if  $B$ is a subfactor with 
$B' \cap M = \mathbb{C}$ then $H(M|B) = \log [M:B]$ by \cite {PP} so that 
$h(A|B) \leq \log [M:B]$. 
\smallskip

St$\o$rmer says that relative entropy can be viewed as a measure of distance between 
two subalgebras, which in the noncommutative case also measures their sizes and relative 
position. 

Here, we give an example, which shows that $h(A | B)$ measures relative 
position and that some small size subalgebra $A_G$ can take the maximal value of $h(A | B),$ 
although  the entropy $h(A | B)$ is increasing in $A.$ 
\vskip 0.5cm

Assume that the finite group $G$ in Theorem 3.1.4 is  abelian. Let $B = uNu^*.$ 
By taking the inner automorphism $Ad u^*,$ we may consider 
$M$  as the crossed product of $B$ by  $G$ so that 
$x \in M$  has a unique expansion $x = \sum_{g \in G} x_gv_g, \ (x_g \in B).$ 
Let $A_G$  be the von Neumann algebra generated by the unitary group $v_G,$ 
(that is, $|G|$ dimensional abelian algebra).  
Then 
$$h(A_G|B) = \log |G| = H(M|B) = \log[M:B].$$ 
\smallskip

In fact, it is clear that 
$h(A|B) \leq H(M|B) = \log[M:B] = \log |G|.$ 
To show the opposite iniquality, 
let $\hat G$ be the character group of $G.$ 
Given $\chi \in \hat G,$ let 
$$p_\chi = \frac 1{|G|} \sum_{g \in G} \chi_g v_g.$$ 
Then $\{p_\chi ; \chi \in \hat G\}$ is a family of 
mutually orthogonal projections in $A_G$ 
with $\sum_{\chi \in \hat G} p_\chi = 1.$ 
Hence 
$${h( A_G \mid B) 
\geq \sum_{\chi \in \hat G} (\tau \eta E_{B}(p_\chi)) 
  = \sum_{\chi \in \hat G} \eta (\frac 1{|G|}) 
  = |\hat G| \frac 1{|G|} \log( {|G|} ) 
 = \log |G|},$$
 and we have that 
$${h( A_G \mid B) =  \log |G|}.$$
\vskip 0.3cm

In the next section, we show that inner conjugacy classes of subfactors $N$ of type II$_1$ factor 
can take the maximum value of $h(N | uNu^*).$ 

\vskip 0.3cm
\subsection{Case of $G = {\mathbb{Z}}_n$ }
Here, we assume that the  group $G$ in 3.1 is a finite cyclic group ${\mathbb{Z}}_n$,  
that is,  
$M$ is the crossed product $N \rtimes_\alpha {\mathbb{Z}}_n$ of a II$_1$-factor $N$ 
by the group generated by an  automorphism  $\alpha$  on $N$ 
such that 
$\alpha^n$ is the identity and $\alpha^i$ is outer for $i = 1, \cdots, n-1.$ 
Such an automorphism  $\alpha$ is called a minimal periodic automorphism (cf. \cite {Co1}). 
\vskip 0.3cm

\subsubsection{Matrix units  for  minimal periodic automorphisms} 
Let $\gamma$ be a  primitive $n$-th root. 
Connes showed in the proof for the characterization of minimal periodic automorphisms 
(\cite [{Cor. 2.7}] {Co1}) that if $\alpha$  is minimal periodic, then there exists 
a set of  matrix units $\{e_{ij}\}_{i,j =1}^n$ in $N$ 
such that 
$$\alpha(e_{ij}) = \gamma^{i-j} e_{ij}, \quad (i,j = 1, \cdots, n).$$

Let $w = \sum_i e_{i+1, i}.$ 
Then $w$ is a unitary in $N$ which satisfies that 
$$w^j = \sum_i e_{i+j, i}, \quad  w^{i*} e_{jj}w^i = e_{j-i,j-i} \quad  \text{and}  
\quad  \alpha(w) = \gamma w.$$ 
\smallskip

The following indicates that 
the inner congugacy class of $N$ can take the $\max_{L} h(N | L),$  
where $L$ is a subfactor of $M$ with $[M : L] = n.$ 
\vskip 0.3cm

\noindent{\bf Theorem 3.2.2.} 
{\it Let $N \subset M$ be the above. Then there exists a  unitary operator 
$u$ in $M$ which satisfies the following properties : }
\smallskip

\noindent(1) $h(N | u N u^*) = H_N(\text{Ad}u) = \log n.$ 
\smallskip

\noindent(2) {\it The conditional expectations $E_N$ and $E_{u N u^*}$ commute. }
\smallskip

\begin{proof} 
Let $w \in N$ be the unitary operator in 3.2.1. Let 
$v \in M$ be a unitary in $M$  implimenting $\alpha,$ that is, 
$\alpha(x) = vxv^*$ for all $x \in N.$ 
We put 
$$u = \frac 1 {\sqrt{n}} \sum_i w^i v^{i-1}.$$ 
Then $u$ is a unitary and 
$\sqrt n E_N(u) = w.$ 
Since  $H_N(u) = \log n$,      
we have  by  Theorem 3.1.4 that 
$$h(N | u N u^*) \leq \log n.$$ 
As a finite partition of the unity, we choose 
$\{p_j : p_j = e_{jj}, j = 1, \cdots, n\}.$ 
Then 
\begin{eqnarray*} 
\lefteqn{h(N | u N u^*)}\\
&&\geq 
\sum_j \tau \eta E_{u N u^*} (p_j) 
- \tau \eta(p_j)  
=  \sum_j \tau \eta E_N( u^* p_j u) \\
&& = \sum_j \tau \eta (\sum_i \alpha^{-i}(\frac{w^{i*}p_jw^i } n)) 
 = \sum_j \tau \eta(\frac{\sum_i p_{j-i}} n)) \\
&&= \log n
\end{eqnarray*}
Hence $h(N | u N u^*) = \log n.$ 
\smallskip 

(2) To show that $E_{u N u^*} E_N = E_N E_{u N u^*},$ 
remark that for all $a \in N,$ 

$$E_{u N u^*}(av^k) = 
\frac 1{n^2} \sum_j ( \sum_{i, l} \gamma^{k^2 + 2ki-(j+l)k-jl} 
w^{j+l-k-i}\alpha^{j+l-k-i}(a) w^{i-l}) v^j.$$ 

Assume that $k \ne 0.$ 
Then $E_N(a v^k) = 0.$ 
On the other hand, 
$$E_N E_{u N u^*}(av^k) 
= \frac 1{n^2}\gamma^{2k}w^{-k}\alpha^{-1}(\sum_i \sum_l 
 \gamma^{2ki-lk} w^{l-i}\alpha^{l-i}(a) w^{i-l})$$ 
and 
$$\sum_i \sum_l  \gamma^{2ki-lk} w^{l-i}\alpha^{l-i}(a) w^{i-l} 
= \sum_j (\sum_i \gamma^{k(i-j)}  w^j \alpha^j(a)w^{*j}) = 0 .$$
Therefore, 
$$E_N E_{u N u^*} (av^k) = 0 = E_{u N u^*} E_N(av^k).$$ 
for all $a \in N$ and $k = 1, \cdots, n-1$. Also for each $a \in N,$ we have that 
\begin{eqnarray*} 
\lefteqn{E_{u N u^*} E_N(a)}\\
&& = \frac 1{n^2} \sum_j( \sum_{i,l}  \gamma^{jl} w^{j+l-i}\alpha^{j+l-i}(a) 
  w^{i-l} ) v^j 
= \frac 1{n^2} \sum_{l} \sum_i w^i \alpha^i(a) w^{*i}\\
&&= \frac 1{n} \sum_i w^i \alpha^i(a) w^{*i}
= E_N(\frac 1{n} \sum_i w^i \alpha^i(a) w^{*i})\\
&&= E_N E_{u N u^*} (a).
\end{eqnarray*}
These show that 
$$E_{u N u^*} E_N(x) = E_N E_{u N u^*} (x) \quad \text{for all} \quad x \in M.$$ 
 \end{proof}
\vskip 0.3cm

\noindent{\bf 3.2.3  \ A continuous family of  subfactors with index 2} 

At the last, in the case of $G = {\mathbb{Z}}_2$, we show 
a  result corresponding  one  in \cite{Ch3} for  maximal abelian subalgebras of 
the type I$_n$ factos $M_n(\mathbb{C})$. 

\vskip 0.3cm

\noindent{\bf Theorem 3.2.3.} 
{\it Let $N$ be a type II$_1$ factor and let $M$ be  the crossed product 
$N \rtimes_\alpha {\mathbb{Z}}_2$ by an outer 
automorphism $\alpha$ with the period 2. 
For the unitary $w \in N$ in  3.2.1,  let 
$$u(\lambda) = \sqrt \lambda w + \sqrt {1-\lambda} v, \quad 
(0 \leq \lambda \leq 1).$$ 
Then $u(\lambda) $ is a unitary in $M$ 
which satisfies the followings : } 
\smallskip 

\noindent
(1) $$h(N \ | \  u(\lambda) \ N \ u(\lambda)^*) = H_N(u) = 
\eta(\lambda) + \eta(1 - \lambda).$$ 
and 
$$\{ h(N \ | \  u(\lambda) \ N \ u(\lambda)^*) : \lambda \in [0, 1] \} = [0, \log 2].$$ 
\smallskip 

\noindent(2) {\it 
\begin{eqnarray*}
\lefteqn{ \quad \quad \ N \quad \quad \quad \quad \quad \subset 
               \quad \quad \quad M }\\
&& \quad \cup  \ \ \ \quad \quad \quad \quad \quad \quad \quad \quad \cup  \\
&& N \cap u(\lambda)Nu(\lambda)^* \  \subset \ u(\lambda)Nu(\lambda)^*
\end{eqnarray*} 
is a commuting square in the sense of \cite{GHJ} if and only if $\lambda = \frac12.$ } 
\smallskip 

\noindent(3) 
$$h(N \ | \  u(\frac 12) \ N \ u(\frac 12)^*) = H_N(Ad u(\frac 12)) =  
\max_{u} h(N \ | uNu^*) = \log 2$$ 
where $u$ is a  unitary in $M.$  
\vskip 0.3cm

\begin{proof} 
(1) It is clear that $H_N(Ad u(\lambda)) = \eta(\lambda) + \eta(1 - \lambda)$. 
Hence by  Theorem 3.1.4, we have 
$$h(N | u(\lambda)Nu(\lambda)^*) \leq \eta(\lambda) 
+ \eta(1 - \lambda).$$ 
We remark that for each $x \in N,$ 
$$E_{u(\lambda)Nu(\lambda)^*}(x) 
= u(\lambda)  E_N( u(\lambda)^* x u(\lambda) )u(\lambda) ^* 
= \lambda w^* x w + (1-\lambda) \alpha(x ). $$
Let $\{e_{ij}\}_{i,j = 1,2} \subset N$ be the set of matrix units for $\alpha$ in 3.2.1.  
Then 
$$E_{u(\lambda)Nu(\lambda)^*}(e_{ii} ) 
= \lambda w^* e_{ii} w + (1-\lambda)    \alpha(e_{ii} ) 
= \lambda  e_{i+1, i+1}  + (1-\lambda)     e_{ii}, \quad (\text{mod}\ 2).$$
Hence, we have that for each $i = 1,2,$ 
$$ \tau \eta (E_{u(\lambda)Nu(\lambda)^*}(e_{ii} ) ) = 
 \frac 12 (\eta(\lambda) + \eta( 1 - \lambda)),$$
so that 
$$h(N | u(\lambda)Nu(\lambda)^*) \geq \eta(\lambda) 
+ \eta(1 - \lambda).$$ 
This implies  that 
$$h(N | u(\lambda)Nu(\lambda)^*) = \eta(\lambda) 
+ \eta(1 - \lambda).$$ 
\vskip 0.3cm

(2) First we remember the following ; the diagram is a  commuting square 
in the sense of \cite{GHJ} means that $E_N E_{uNu^*} =  E_{uNu^*} E_N.$ 

Let $x \in N.$ 
Since $\alpha(w) = -w^*,$ we have that 
$$ E_N E_{u(\lambda) Nu(\lambda)^*}(x)  
= \lambda^2 x + 2 \lambda (1- \lambda) w \alpha(x) w^* 
  + (1 - \lambda)^2 x$$
and 
\begin{eqnarray*} 
\lefteqn{E_{uNu^*}E_{N}(x)} \\
&& = \lambda^2 x + 2 \lambda (1- \lambda) w \alpha(x) w^* 
  + (1 - \lambda)^2 x + \sqrt{\lambda(1-\lambda)}(2\lambda -1) 
  xwv
\end{eqnarray*}
Hence $E_N E_{u(\lambda) Nu(\lambda)^*}(x)  = 
E_{uNu^*}E_{N}(x)$ for all $x \in N$ if and only if 
$\lambda = 1/2.$ 
Similarly, for all $x \in N,$  
$$ E_N E_{u(1/2) Nu(1/2)^*}(xv)  = w \alpha(x_1) + x_1w^* + x_1 \alpha(w) + \alpha(w^*x_1)  = 0 $$
and 
$$E_{u(1/2)Nu(1/2)^*}E_{N}(xv) = 0.$$
These imply the conclusion. 
\smallskip 

\noindent(3) 
Since $N$ is a subfactor of $M$ with $[M : N] = 2,$  we have that  
$$h(N | uNu^*) \leq H(N | uNu^*)  \leq H(M |N) = \log 2$$ 
so that 
$$h(N \ | \  u(\frac 12) \ N \ u(\frac 12)^*) =  \log 2 =  
\max_{u} h(N \ | uNu^*).$$
\end{proof}

\thanks{The author was supported in part by JSPS Grant No.20540209.}
\end{document}